\documentclass[10pt]{article}

\usepackage{amsmath}
\usepackage{amsfonts}
\usepackage{amssymb}

\let\ds=\displaystyle
\def\R{\mathbb{R}_q}

\def\N{\mathbb{N}}
\def\Z{\mathbb{Z}}

\def\F{\mathcal{F}_{q,\v}}
\def\I{\infty}
\def\ll{\mathcal{L}_{q,1,\v}}
\def\cc{\mathcal{C}_{q,0}}
\def\cb{\mathcal{C}_{q,b}}
\def\lL{\mathcal{L}_{q,2,\v}}
\def\lp{\mathcal{L}_{q,p,\v}}

\def\lr{\mathcal{L}_{q,r,\v}}
\def\Lp{\mathcal{L}_{q,p',\v}}

\def\A{\mathcal{A}_{q,\v}}

\def\v{\nu}

\newtheorem{theore}{Theorem}
\newtheorem{coro}{Corollary}
\newtheorem{definitio}{Definition}

\newtheorem{lemm}{Lemma}
\newtheorem{propositio}{Proposition}
\newtheorem{remar}{Remark}
\newenvironment{proo}[1][Proof]{\noindent\textbf{#1.} }{\ \rule{0.5em}{0.5em}}

\begin{document}

\title{\bf $q$-Macdonald function as a Variation Diminishing $*_q$-kernel}

\author{
Lazhar Dhaouadi\thanks{Institut Pr\'eparatoire aux Etudes d'Ing\'enieur de Bizerte, 7021 Zarzouna,Tunisia. E-mail : lazhardhaouadi@yahoo.fr},\quad
Islem Saidani\thanks{Faculty of Sciences of Bizerte, 7021 Zarzouna,Tunisia. E-mail : saidani.islem1@gmail.com}\quad and\quad
Hedi  Elmonser\thanks{College of Science and Humanities at Howtat Sudair, Majmaah University, Saudi Arabia. E-mail : monseur2004@yahoo.fr}
}
\date{ }
\maketitle

\begin{abstract}

In this paper we study the variation diminishing kernel as a part of the $q$-calculus. We introduce the $q$-Macdonald function a newborne in the family of the $q$-special functions which play a central role in this study.

\vspace{5mm}

\vspace{3mm}
\end{abstract}

\section{Introduction}

We denote by $V[a_1,\ldots,a_n]$ the number of variation of sign of the sequence $a_1,\ldots,a_n$ of real numbers. If $f$ is real function defined on an interval $(a,b)$, we define $V[f]$ the variation of $f$ on $(a,b)$, as
$$
V[f]=\sup V[f(a_1),\ldots,f(a_n)]
$$
where the supremum is taken over all finite lists $a<a_1,\ldots,a_n<b$. It is possible of course that $V[f]=\I$.\\

Finally, a real-valued kernel $T(x,y)$ define on $(a,b)\times (a,b)$ is said to be variation diminishing kernel if
$$
V[Tf(x)]\leq V[f(x)]
$$
for all real-valued integrable function $f$, where
$$
Tf(x)=\int_a^b T(x,y)f(y)dy.
$$
Investigations and efforts have been provided by several authors to study this vast subject. See Karlin \cite{K2} for further information about variation diminishing transformations.\\

A real-valued function $\varphi\in L^1(a,b)$ is said to be a variation diminishing convolution kernel if  $T(x,y)=\varphi(x-y)$ is a variation diminishing kernel. In this cases $Tf=\varphi*f$.\\

The most importent think about the variation diminishing convolution kernel is their connection with the Laguerre-Polya-Schur class. This connection  was firstly proved by I. J. Schoenberg  \cite{S2,S1} for the convolution product associated to the Fourier transform and later  by I. I. Hirshman \cite{H} for the Hankel transform.\\

If we can qualify a special function as being important when it appears in mathematical and physical applications, then the modified Bessel function of the third kind is a quite important one. For instance, the modified Bessel function  $K_\nu(x)$ is one of the solutions to the modified Bessel differential equation. Also known as modified Hankel function or Macdonald function. Inspired by information about the  Macdonald function which was found in the book by Watson \cite{W} an analogue of this function is given as a part of the $q$-calculus. However, our interest motivated by the various properties of the $q$-Macdonald  function $K_{\nu,q}(x)$ when we  study it as a variation diminishing $q$-convolution kernel. This approche is original and not be fond in the literature of classical harmonic analysis.\\

This paper is organized as follows: Section $2$ is devoted to an overview of the $q$-Bessel Fourier transform. In section $3$, we introduce the
concept of variation diminishing $q$-convolution kernel. Also we define the $q$-Macdonald  function and we study some of it's properties. In particular we prove that the $q$-Macdonald  function is a variation diminishing $q$-convolution kernel and as a consequence it's a non negative function on $\R$. In section 4 we deduce some result about composite variation diminishing $q$-convolution kernel. We prove that the $q$-Gauss kernel $e^c_{\v,q}$ is a variation diminishing $*_q$-kernel. In section 5 we discuss the asymptotic expansion at infinity of the variation diminishing $q$-convolution kernel.

\section{Preliminaries about the $q$-Bessel Fourier transform}
Throughout this paper we adopt the standard conventional notations of \cite{G}.   Let $0<q<1,\quad\v>-1$ and consider
$$
\mathbb{R}_q=\{q^n,\quad n\in\Z\},\quad \R^+=\{q^n,\quad n\in\N\},\quad\mathbb{R}_q^-=\{q^n,\quad n\in\Z\setminus\N\}.
$$
 For any complex number $a$
$$
(a;q)_0=1,\quad (a;q)_n=\prod_{i=0}^{n-1}(1-aq^{i}), \quad n\in\N^*.
$$
The $q$-derivative of a function $f$ is defined  for $x\neq 0$ by
$$
D_qf(x)=\frac{f(x)-f(qx)}{(1-q)x}
$$
and we have
$$
D_q(fg)=(D_qf)g+(\Lambda_qf)D_qg,
$$
where $\Lambda_q$ is the $q$-shift operator defined by $\Lambda_qf(x)=f(qx)$.\\

The Jackson's $q$-integrals are defined by \cite{J}
$$
\int_0^1f(x)d_qx=(1-q)\sum_{n=0}^\infty q^nf(q^n)=(1-q)\sum_{x\in\R}xf(x),
$$

$$
\int_0^\infty f(x)d_qx=(1-q)\sum_{n=-\infty}^\infty q^nf(q^n)=(1-q)\sum_{x\in\mathbb{R}_q} xf(x),
$$

$$
\int_1^\I f(x)d_qx=(1-q)\ds\sum_{n\in\Z\setminus\N}q^nf(aq^n)=(1-q)\ds\sum_{x\in\mathbb{R}_q^-}xf(x).
$$
Also we have the following identity
$$
\int_0^\I f(ax)g(x)d_qx=\frac{1}{a}\int_0^\I f(x)g(x/a)d_qx,\quad\forall a\in\R,
$$
provide $\ds\int_0^\I f(ax)g(x)d_qx$ exists.\\

The  normalized $q$-Bessel function of Hahn-Exton is defined as follows \cite[p.43]{D}
\begin{equation*}
j_{\v}(x,q^{2})=\sum_{n=0}^{\infty }(-1)^{n}\frac{q^{n(n+1)}}{
(q^{2\v+2},q^{2})_{n}(q^{2},q^{2})_{n}}x^{2n}.
\end{equation*}
It satisfies the following estimate \cite[p.44]{D1}
\begin{equation*}
|j_\v(q^{n},q^{2})|\leq \frac{(-q^{2};q^{2})_{\infty
}(-q^{2\v+2};q^{2})_{\infty }}{(q^{2\v+2};q^{2})_{\infty }}\left\{
\begin{array}{c}
1\quad \quad \quad \quad \quad \text{if}\quad n\geq 0 \\
q^{n^{2}-(2\v+1)n}\quad \text{if}\quad n<0
\end{array}
\right..
\end{equation*}
The function $x\mapsto j_\v(tx,q^2)$ is a solution of the following $q$-differential equation \cite[p.43]{D}
\begin{equation}\label{e7}
\Delta_{q,\v}f(x)=-t^2f(x),
\end{equation}
where $\Delta _{q,\v}$ is the $q$-Bessel operator
\begin{equation}\label{e10}
\Delta _{q,\v}f(x)=\frac{1}{x^{2}}\Big[ f(q^{-1}x)-(1+q^{2\v})f(x)+q^{2\v}f(qx)\Big].
\end{equation}
The $q$-Wronskian was introduced in \cite[p.58]{D2} as follows
$$
w_x(f,g)=q^{-1}(1-q)^2\Big[\Lambda_q^{-1}D_qf(x)g(x)-\Lambda_q^{-1}D_qf(x)g(x)\Big],
$$
and we have
\begin{equation}\label{e18}
D_q\Big[y\mapsto
y^{2v+1}w_y(f,g)\Big](x)=\Big[\Delta_{q,\v}f(x)g(x)-f(x)\Delta_{q,\v}g(x)\Big]x^{2\v+1}.
\end{equation}

Also the  normalized $q$-Bessel function $j_{\nu }(.,q^{2})$ satisfies the following orthogonality relation \cite[p.43]{D}:
\begin{equation}\label{c}
c_{q,\nu }^{2}{\int_{0}^{+\infty }}j_{\nu }(xt,q^{2})j_{\nu
}(yt,q^{2})t^{2\nu +1}d_{q}t=\delta _{q}(x,y),\text{ \ \ \ }\forall x,y\in
\mathbb{R}_{q}^{+}
\end{equation}%
where
$$
c_{q,\v}=\frac{1}{1-q}\frac{(q^{2\v+2};q^2)_\infty}{(q^{2};q^2)_\infty},
$$
and
\begin{equation*}
\delta _{q}(x,y)=\left\{
\begin{array}{c}
0~~\ ~~~~~~~~~\ \ ~~~~~\text{si }x\neq y\  \\
\frac{1}{(1-q)x^{2(\nu +1)}}~~~~~\text{si }x=y
\end{array}%
\right. .
\end{equation*}
For the particulars cases $x=q^{n}~$ and $y=q^{m}$ in (\ref{c}) we have
\begin{equation*}
c_{q,\nu }^{2}\int_{0}^{\infty }j_{\nu }(tq^{n},q^{2})\times j_{\nu}(tq^{m},q^{2})t^{2\nu +1}d_{q}t=\frac{q^{-2n(\nu+1)}}{1-q}\delta_{nm}.
\end{equation*}
When $n \rightarrow +\infty ~$ and by the dominate convergence Theorem we obtain
\begin{equation}\label{or}
\int_{0}^{\infty }j_{\nu }(q^{m}t,q^{2})t^{2\nu +1}d_{q}t=0.
\end{equation}
The $q$-Bessel Fourier transform $\mathcal{F}_{q,\v}$ is defined by \cite[p.44]{D}
\begin{equation*}
\mathcal{F}_{q,\v}f(x)=c_{q,\v}\int_{0}^{\infty
}f(t)j_{\v}(xt,q^{2})t^{2\v+1}d_{q}t,\quad\forall x\in \mathbb{R}_{q}.
\end{equation*}
We introduce some functional spaces:\\
\begin{quote}
$\bullet$ The space $\lp$ , $1\leq p<\I$ denote the sets of real functions on $\R$ for which
$$
\|f\|_{q,p,\v}=\left[\int_0^\I |f(x)|^px^{2\v+1}d_qx\right]^{1/p}
$$
is finite. The inner product on the Hilbert space $\lL$ is defined by:
$$
\langle f,g\rangle=\int_0^\I f(x)g(x)x^{2\v+1}d_qx.
$$

$\bullet$ The spaces $\cc$ and $\cb$ denote the set of functions defined on $\R$ and $\ds\lim_{n\to\I}f(q^n)$ exists, which are respectively
vanishing at infinity and bounded. These spaces are equipped with the topology of uniform convergence.\\

$\bullet$ The $q$-Wiener algebra $\A=\Big\{f\in\ll,\quad \F(f)\in\ll\Big\}$ is a subspace of $\cb$ and
\begin{equation}\label{e8}
\overline{\A}=\cb.
\end{equation}
\end{quote}

Given a function $f\in\mathcal L_{q,1,\v}$ then  ${\mathcal F}_{q,\v}f\in{\mathcal C}_{q,0}$ and \cite[p.44]{D}
\begin{equation}\label{e11}
\|\mathcal F_{q,\v}f\|_{q,\infty}\leq B_{q,\v}\|f\|_{q,1,\v}
\end{equation}
where
$$B_{q,\v}=\frac{1}{1-q}\frac{(-q^2;q^2)_\infty(-q^{2\v+2};q^2)_\infty}{(q^2;q^2)_\infty}.$$
Given $f\in\ll$ then we have the following inversion formula \cite[p.45]{D}
\begin{equation}\label{e1}
\F^2(f)(x)=f(x),\quad\forall x\in\R.
\end{equation}
The $q$-Bessel Fourier transform $ \F:\lL\rightarrow\lL$ defines an isomorphism and we have\cite[p.45]{D}
\begin{equation*}
\F^2(f)=f,\quad\|\F(f)\|_{q,2,\v}=\|f\|_{q,2,\v}.
\end{equation*}
Let $\v>0$ and suppose that the following integrals are finite \cite[p.60]{D2}
$$
\int_0^\I f(qx)g(x)x^{2\v-1}d_qx,\quad \int_0^\I
f(x)g(x)x^{2\v-1}d_qx,\quad \int_0^\I f(x/q)g(x)x^{2\v-1}d_qx,
$$
then
\begin{equation}\label{e22}
\langle\Delta_{q,\v}f,g\rangle=\langle f,\Delta_{q,\v}g\rangle.
\end{equation}
In particular if $f\in\ll$ and $\ds\lim_{x\to 0}f(x)$ exists then $\Delta_{q,\v}f\in\ll$ and we have
\begin{equation}\label{e23}
\F(\Delta_{q,\v}f)(x)=-x^2\F(\Delta_{q,\v}f)(x),\quad\forall x\in\R.
\end{equation}
Note that (\ref{e22}) hold true if we moreover require that \cite[p.60]{D2}
\begin{equation}\label{e24}
D_qf(x)=O(x^{-\v})\quad\text{and}\quad D_qg(x)=O(x^{-\v}),
\end{equation}
as $x\downarrow 0$.\\

The $q$-Bessel translation operator is given by \cite[p.47]{D}
\begin{equation*}
T_{q,x}^{\v}f(y)=c_{q,\v}\int_{0}^{\infty }\mathcal{F}%
_{q,\v}f(t)j_{\v}(yt,q^{2})j_{\v}(xt,q^{2})t^{2\v+1}d_{q}t.
\end{equation*}
This operator can be written as follows
\begin{equation*}
T_{q,x}^{\v}f(y)=\int_{0}^{\infty }f(z)D_{q,\v}(x,y,z)z^{2\v+1}d_{q}z,
\end{equation*}
with  kernel
\begin{equation*}
D_{q,\v}(x,y,z)=c_{q,\v}^{2}\int_{0}^{\infty
}j_{\v}(xs,q^{2}j_{\v}(ys,q^{2}j_{\v}(zs,q^{2})s^{2\v+1}d_{q}s.
\end{equation*}
The $q$-convolution product is given by \cite[p.49]{D}
\begin{equation*}
f\ast _{q}g(x)=c_{q,\v}\int_{0}^{\infty
}T_{q,x}^{\v}f(y)g(y)y^{2\v+1}d_{q}y.
\end{equation*}
Given two functions $f,g\in\ll$ then \cite[p.49]{D}
\begin{equation}\label{e5}
f*_qg\in\ll,
\end{equation}
and
\begin{equation}\label{e6}
\F(f*_qg)=\F(f)\times\F(g).
\end{equation}
Let $1\leq p,p',r$ such that $\ds{1\over p}+{1\over p'}-1={1\over r}.$ If $f\in\lp$ and $g\in\Lp$ then
\begin{equation*}
f*_qg\in\lr
\end{equation*}
and
\begin{equation*}
\|f*_qg\|_{q,r,\v}\leq B_{q,p,\v}B_{q,p',\v}B_{q,r',v}\|f\|_{q,p,\v}\|g\|_{q,p',\v}
\end{equation*}
where
$$
\frac{1}{r}+\frac{1}{r'}=1,\quad
B_{p,q,\v}=B_{q,\v}^{(\frac{2}{p}-1)}.
$$
In particular if $f\in\ll$ and $g\in \cb$ then $f*_qg\in\cb$ and we have
\begin{equation}\label{e12}
\|f*_qg\|_{q,\I}\leq \|f\|_{q,1,\v}\|g\|_{q,\I}.
\end{equation}
Similarly if $f\in\ll$ and $g\in\ll$ then
\begin{equation*}
\|f*_qg\|_{q,1,\v}\leq \|f\|_{q,1,\v}\|g\|_{q,1,\v}.
\end{equation*}
The $q^2$-exponential  function \cite[p.]{} is defined on $\R$ by
\begin{equation*}
e(-x^2,q^2)=\frac{1}{(-x^2;q^2)_\I}.
\end{equation*}
The $q$-Gauss kernel is defined by \cite[p.6]{D1}
\begin{equation}\label{e17}
e^c_{\v,q}(x)=\F\Big[t\mapsto e(-c^2t^2,q^2)\Big](x)=\frac{(-q^{2\v+2}c,-q^{-2\v}/c;q^2)_\I}{(-c,-q^2/c;q^2)_\I}~e\left(-\frac{x^2}{q^{2\v}c},q^2\right),\quad c>0.
\end{equation}
Given $f\in\ll\cap\lp$ where $1\leq p<\infty$, then \cite[p.8]{D1}
\begin{equation}\label{e9}
\lim_{c\to  0}\|f-f*_qe^c_{\v,q}\|_{q,p,\v}=0.
\end{equation}

\bigskip

\section{Elementary kernel}

A real function $f$ defined on $\mathbb{R}_q$ is said to have at least $n$ changes of sign if there exists numbers
$$
0<t_0<t_1<\ldots<t_n,\quad t_i\in\mathbb{R}_q
$$
such that
$$
f(t_i)f(t_{i-1})<0,\quad i=1,\ldots,n.
$$
The function $f$ has exactly $n$ changes of sign if it has at least $n$ changes of sign and does not have at least $n+1$ changes of sign. The number
of changes of sign of $f$ is denoted by $V[f]$ which has one of the values $0,1,\ldots$ or $+\infty$.

\begin{definitio}
A function $K\in\mathcal{L}_{q,1,\v}$ is said to be a variation diminishing $*_q$-kernel if for every $f\in\mathcal{C}_{q,b}$ we
have $ V[K*_qf]\leq V[f].$
\end{definitio}

\begin{lemm}\label{lem1}
If $f$ is a function defined on $\mathbb{R}_q$ such that either
$$
\lim_{x\rightarrow 0^+}f(x)=0,\quad\text{or}\quad\lim_{x\rightarrow
+\infty}f(x)=0
$$
holds, then
\begin{equation}\label{e}
V\Big[D_qf\Big]\geq V[f].
\end{equation}
\end{lemm}

\begin{proo}
\begin{quote}
. If $f$ change sign at $x$ and $f(x)>0$ then $D_qf(x)>0$.

. If $f$ change sign at $x$ and $f(x)<0$ then $D_qf(x)<0$.\\
\end{quote}
This means that for two consecutive change of sign of $f$ there is one change of sign of $D_qf$. So if $V(f)=+\infty$, we have $V\Big[D_qf\Big]=+\I$ and  then inequality (\ref{e}) is true.\\

Let's look at the cases when $V(f)=n\geq 2$. Let $\{q_1,\ldots q_n\}$ the number when $f$ change sign. Hence at least $D_qf$ has $(n-1)$ change of sign.\\

{\bf case i.} $\ds\lim_{x\to \I} f(x)=0$:\\

 suppose $f>0$ on  $[q_n,+\I)$. The function $f$ should be decrease toward $0$ when $x\to+\I$. Then there exists
at least a number $q_{n+1}>q_n$ such that $D_qf(q_{n+1})<0$. As $D_qf(q_{n})>0$ we deduce inequality (\ref{e}). Similar argument if $f<0$ on  $[q_n,+\I)$.\\

{\bf case ii.} $\ds\lim_{x\to 0^+} f(x)=0$:\\

suppose $f>0$ on  $(0,q_1)$. The function $f$ should be decrease toward $0$ when $x\downarrow  0^+$. Then there exists
at least a number $q_0<q_1$ such that $D_qf(q_0)<0$. As $D_qf(q_1)>0$ we deduce inequality (\ref{e}). Similar argument if $f<0$ on  $(0,q_1)$.
\end{proo}

\begin{coro}\label{c0}
Let $h$ a function defined on $\mathbb{R}_q$ and $\Omega$ be a positive  function such that either
$$
\lim_{x\rightarrow 0^+}\Omega(x)h(x)=0
\quad\text{or}\quad\lim_{x\rightarrow +\infty}\Omega(x)h(x)=0
$$
holds, then $ V\Big[D_q(\Omega h)\Big]\geq V[h].$
\end{coro}

\begin{remar}\label{r0}Let $(h_n)_n$ a sequence of function defined on $\R$ which converge pointwise  to $h$
$$
\lim_{n\to\I}|h_n(x)-h(x)|=0,\quad\forall x\in\R.
$$
Let $\{q_1,\ldots q_r\}$ a given number of $\R$. For any $\epsilon>0$ there is $n_0\in\N$ such that
$$
|h_n(q_i)-h(q_i)|<\epsilon,\quad i=1\ldots r,\quad\forall n\geq n_0.
$$
If we take $\epsilon<\min(|h(q_i)|,\quad i=1\ldots r)$ then
$$
{\rm sign}~h_n(q_i)={\rm sign}~h(q_i),\quad i=1\ldots r.
$$
There fore
$$
V[h_n]= V[h],\quad\forall n\geq n_0\Rightarrow \lim_{n\to\I}V[h_n]= V[h],
$$
provide $V[h]$ is finite, if note $\ds\lim_{n\to\I}V[h_n]= V[h].=+\I$.
\end{remar}
In analogy to the integral representation of the Macdonald function in \cite[p. 434]{W} we give a $q$-version of this function.

\begin{definitio} The $q$-Macdonald function is defined by
$$
K_{\v,q}(x)=c_{q,\v}\int_0^\I\left[1+t^2\right]^{-1}j_\v(tx,q^2)t^{2\v+1}d_qt.
$$
\end{definitio}

For $a\in\R$, let us set
$$
K^a_{\v,q}(x)=c_{q,\v}\int_0^\I\left[1+\frac{t^2}{a^2}\right]^{-1}j_\v(tx,q^2)t^{2\v+1}d_qt.
$$
In particular $K_{\v,q}^a(x)=a^{2(\v+1)}K_{\v,q}(ax).$

\begin{theore}\label{t1}
The function $K^a_{\v,q}\in\ll$ and we have
\begin{equation}\label{e2}
\mathcal{F}_{q,\v}(K^a_{\v,q})(x)=\left[1+\frac{x^2}{a^2}\right]^{-1},\quad\forall
x\in\mathbb{R}_q.
\end{equation}
In addition, the function $K_{\v,q}^a$ satisfies the following $q$-difference
equation
\begin{equation}\label{e4}
\left[1-\frac{\Delta_{q,\v}}{a^2}\right]K_{\v,q}^a(x)=0.
\end{equation}
\end{theore}

\begin{proo}First, we will prove that
\begin{equation}\label{e0}
\int_0^1|K_{\v,q}^a(x)|x^{2\nu+1}d_qx<\I.
\end{equation}
{\bf case i.} $\v\geq-1/2$: we have
\begin{eqnarray*}
x^{2\v+2}|K_{\v,q}^a(x)|&\leq&c_{q,\v}\int_0^\I
\left[1+\frac{t^2}{a^2x^2}\right]^{-1}|j_\v(t,q^2)|t^{2\v+1}d_qt\\
&\leq&\left[c_{q,\v}\int_0^\I
\sqrt{\left(\frac{t}{ax}\right)}\left[1+\frac{t^2}{a^2x^2}\right]^{-1}|j_\v(t,q^2)|t^{2\v+1/2}d_qt\right]\sqrt{ax}\\
&\leq&\left[c_{q,\v}\int_0^\I |j_\v(t,q^2)|t^{2\v+1/2}d_qt\right]\sqrt{ax}.
\end{eqnarray*}
So
\begin{equation}\label{e21}
x^{2\v+2}|K_{\v,q}^a(x)|=O(\sqrt{x})
\end{equation}
and $\ds\sum_{x\in\R^+}\sqrt{x}<\I$  which prove (\ref{e0}). \\

{\bf case ii.} $-1<\v<-1/2$:
$$
\int_0^\I\left[1+\frac{t^2}{a^2}\right]^{-1}t^{2\v+1}d_qt<\I.
$$
By the dominate convergence theorem $\ds\lim_{x\to 0}K_{\v,q}^a(x)$ exists. Then (\ref{e0}) hold true.\\

Second, proving
\begin{equation}\label{e20}
\int_1^{\I}|K_{\v,q}^a(x)|x^{2\nu+1}d_qx<\I,\quad\v>-1.
\end{equation}
We have
\begin{eqnarray*}
x^2K_{\v,q}^a(x)&=&c_{q,\v}\int_0^\I
\left[1+\frac{t^2}{a^2}\right]^{-1}x^2j_\v(xt,q^2)t^{2\v+1}d_qt\\
&=&-c_{q,\v}\int_0^\I
\left[1+\frac{t^2}{a^2}\right]^{-1}\Delta_{q,\v}j_\v(xt,q^2)t^{2\v+1}d_qt\\
&=&-c_{q,\v}\int_0^\I
\Delta_{q,\v}\left[1+\frac{t^2}{a^2}\right]^{-1}j_\v(xt,q^2)t^{2\v+1}d_qt,\quad (*)\\
&=&-c_{q,\v}\int_0^\I
u(t)j_\v(xt,q^2)t^{2\v+1}d_qt\\
\end{eqnarray*}
where $u(t)=\Delta_{q,\v}\left[1+\frac{t^2}{a^2}\right]^{-1}$ a bounded function
on $\R$. To justify $(*)$ we use (\ref{e22}) if $\v>0$ and (\ref{e24}) if $-1<\v\leq 0$.
\begin{eqnarray*}
x^2\Big[x^{2\v+2}|K_{\v,q}^a(x)|\Big]&\leq&c_{q,\v}\int_0^\I
|u(t/x)||j_\v(t,q^2)|t^{2\v+1}d_qt\\
&\leq&c_{q,\v}\|u\|_{q,\I}\int_0^\I
|j_\v(t,q^2)|t^{2\v+1}d_qt.
\end{eqnarray*}
Hence $x^{2\v+2}|K_{\v,q}^a(x)|=O(1/x^2)$ which prove (\ref{e20}).\\

To prove the second result (\ref{e4}) we use identity (\ref{or})
\begin{eqnarray*}
\left[1-\frac{\Delta_{q,\v}}{a^2}\right]K^a_{\v,q}(x)&=&c_{q,\v}\int_0^\I \left[1+\frac{t^2}{a^2}\right]%
^{-1}
\left[1-\frac{\Delta_{q,\v}}{a^2}\right]j_\v(tx,q^2)t^{2\v+1}d_qt\\
&=&c_{q,\v}\int_0^\I j_\v(tx,q^2)t^{2\v+1}d_qt=0.
\end{eqnarray*}
Note that
\begin{equation}\label{e16}
\left[1-\frac{\Delta_{q,\v}}{a^2}\right]j_\v(tx,q^2)=\left[1+\frac{t^2}{a^2}
\right]j_\v(tx,q^2).
\end{equation}
This achieves the proof.
\end{proo}

\begin{coro}\label{c1}
If $f\in\mathcal{L}_{q,1,\v}$ and $ h=K^a_{\v,q}*_qf$ then
$$
\left[1-\frac{\Delta_{q,\v}}{a^2}\right]h(x)=f(x).
$$
\end{coro}

\begin{proo}
By (\ref{e5}) and (\ref{e6}) we see that $h\in\ll$ and we have
$$
\mathcal{F}_{q,\v}(h)(t)=\mathcal{F}_{q,\v}(f)(t)\left[1+\frac{t^2}{a^2}\right]%
^{-1}.
$$
By (\ref{e1}), we have
$$
h(x)=c_{q,\v}\int_0^\I \mathcal{F}_{q,\v}(f)(t)\left[1+\frac{t^2}{a^2}\right]%
^{-1} j_\v(tx,q^2)t^{2\v+1}d_qt,
$$
and so
$$
\left[1-\frac{\Delta_{q,\v}}{a^2}\right]h(x) =c_{q,\v}\int_0^\I \mathcal{F}%
_{q,\v}(f)(t)\left[1+\frac{t^2}{a^2}\right]^{-1} \left[1-\frac{\Delta_{q,\v}}{%
a^2}\right]j_\v(tx,q^2)t^{2\v+1}d_qt.
$$
Using (\ref{e16}) and the inversion formula (\ref{e1}) we obtain
$$
\left[1-\frac{\Delta_{q,\v}}{a^2}\right]h(x)=c_{q,\v}\int_0^\I \mathcal{F}%
_{q,\v}(f)(t)j_\v(tx,q^2)t^{2\v+1}d_qt=f(x).
$$
This proves the result.
\end{proo}

\begin{definitio}The modified $q$-Bessel function is defined as follows
$$I_{\v,q}(x)=j_\v(ix,q^2),\quad i^2=-1.$$
\end{definitio}

\begin{remar}\label{re}
The function $I^a_{\v,q}:x\mapsto I_{\v,q}(ax)$ satisfies $\ds\left[1-\frac{\Delta_{q,\v}}{a^2}\right]I^a_{\v,q}(x)=0$ and $I^a_{\v,q}(x)>0;\quad\forall x\in\R.$
\end{remar}

\begin{propositio}\label{p1}
Let $h$ a given function defined on $\R$ then we have
$$
\left[1-\frac{\Delta_{q,\v}}{a^2}\right]h(x) =-\frac{q^{2\v-1}(1-q)^2}{a^2x^{2\v+1}I^a_{\v,q}(x)%
}\Lambda_q^{-1}D_q\left[x^{2\v+1}I^a_{\v,q}(x)I^a_{\v,q}(qx)D_q\left[\frac{h(x)}{I^a_{\v,q}(x)}
\right]\right],
$$
\end{propositio}

\begin{proo}
In fact
\begin{eqnarray*}
&&-\frac{q^{2\v-1}(1-q)^2}{a^{2}x^{2\v+1}I^a_{\v,q}(x)}\Lambda _{q}^{-1}D_{q}\left[
x^{2\v+1}I^a_{\v,q}(x)I^a_{\v,q}(qx)D_{q}\left[ \frac{h(x)}{I^a_{\v,q}(x)}\right] \right] \\
&=&-\frac{q^{2\v-1}(1-q)}{a^{2}x^{2\v+1}I^a_{\v,q}(x)}\Lambda _{q}^{-1}D_{q}\left[
x^{2\v+1}I^a_{\v,q}(x)I^a_{\v,q}(qx)\left[ \frac{\frac{h(x)}{I^a_{\v,q}(x)}-\frac{h(qx)}{
I^a_{\v,q}(qx)}}{x}\right] \right] \\
&=&-\frac{1}{a^{2}x^{2}}h(q^{-1}x)-\frac{q^{2\v}}{a^{2}x^{2}}h(qx)+\frac{1}{
a^{2}x^{2}}\left[\frac{I^a_{\v,q}(q^{-1}x)+q^{2\v}I^a_{\v,q}(qx)}{I^a_{\v,q}(x)}\right]
h(x) \\
&=&-\frac{1}{a^{2}x^{2}}h(q^{-1}x)-\frac{q^{2\v}}{a^{2}x^{2}}h(qx)+\frac{1}{%
a^{2}x^{2}}\left[\frac{(1+q^{2\v})I^a_{\v,q}(x)+a^{2}x^{2}I^a_{\v,q}(x)}{I^a_{\v,q}(x)}\right] h(x) \\
&=&\left[ 1-\frac{\Delta _{q,\v}}{a^{2}}\right] h(x).
\end{eqnarray*}
This proves the result.
\end{proo}

\begin{theore}\label{th}
Let $\v>-1/2$. The function $K^a_{\v,q}$ is a variation diminishing $*_q$-kernel.
\end{theore}
\begin{proo}
By Corollary \ref{c1}, if  $f\in\mathcal{L}_{q,1,\v}$ and $h=K_{\v,q}^a*_qf$ then
$$
\left[ 1-\frac{\Delta _{q,\v}}{a^{2}}\right] h(x)=f(x).
$$
{\bf Cases i.} let $f\in\A$, then by the use of Proposition \ref{p1}, Corollary \ref{c0}  and the fact that
$$
\lim_{x\rightarrow\I}1/I^a_{\v,q}(x)=0,\quad \lim_{x\rightarrow
0^+}x^{2\v+1}I^a_{\v,q}(x)I^a_{\v,q}(qx)=0,
$$
we see that $V[f]\geq V[h]$.\\

{\bf Cases ii.} let $f\in\cb$ then  by (\ref{e8}) there exists a sequence of functions
$f_n\in\A$ such that
$$
\lim_{n\to\I}f_n(x)=f(x)\ \text{in}\ \mathcal{C}_{q,b}\ ,
$$
that's mean
$$
\lim_{n\to\I}\|f_n-f\|_{q,\I}=0.
$$
Let $h_n=K^a_{\v,q}*_qf_n$. Using the first case we have, for each $n$,
$$
V[h_n]\leq V[f_n]\leq V[f].
$$
On the other hand, we have by (\ref{e12})
$$
\|h_n-h\|_{q,\I}=\|K_{\v,q}^a*_qf_n-K_{\v,q}^a*_qf\|_{q,\I}\leq\|K_\v^a*-K_{\v,q}^a\|_{q,1,\v}\|f_n-f\|_{q,\I}.
$$
By passage to the limit and Remark \ref{r0}, we obtain $V[h]\leq V[f]$.
\end{proo}

\begin{lemm}\label{lem3}
If $g$ is a variation diminishing $*_q$-kernel then either
$g(x)\geq 0$ or $g(x)\leq 0$ for all $x\in\R$.
\end{lemm}

\begin{proo}By (\ref{e9}) we get $\ds\lim_{c\to 0}\|g-g*_qe^c_{\v,q}\|_{q,1,\v}=0.$ Since $g$ is variation diminishing $*_q$-kernel and $e^c_{\v,q}\in\cb$ we must have $V[g*_qh_c]\leq V[e^c_{\v,q}]=0$. By the definition of the $q$-Jackson integral we have
$$
\|g-g*_qe^c_{\v,q}\|_{q,1,\v}=(1-q)\sum_{x\in\R}x^{2\v+2}|g(x)-g*_qe^c_{\v,q}(x)|.
$$
Then for a given $x\in\R$ we get $\ds\lim_{c\to 0}|g(x)-g*_qe^c_{\v,q}(x)|=0$. Using Remark \ref{r0} we get
$$
V[g]=\lim_{c\to 0}V[g*_qe^c_{\v,q}]=0.
$$
This prove the result.
\end{proo}

\begin{propositio}\label{rem3}Let $\v>-1/2$. We have
$$
K^a_{\v,q}(x)>0,\quad\forall x\in\R.
$$
\end{propositio}

\begin{proo}The function $K^a_{\v,q}$ is a variation diminishing $*_q$-kernel
and
$$
\F(K^a_{\v,q})(0)=c_{q,\v}\int_0^\I
K^a_{\v,q}(t)t^{2\v+1}d_qt=\left[1+\frac{x^2}{a^2}\right]_{x=0}=1.
$$
Then $K^a_{\v,q}(x)\geq 0$ for all $x\in\R$. Using (\ref{e10}) we see that
$$
K^a_{\v,q}(x/q)+q^{2\v}K^a_{\v,q}(qx)=[1+q^{2\v}+(ax)^2]K^a_{\v,q}(x).
$$
Now if there exists $x\in\R$ such that $K^a_{\v,q}(x)=0$ then
$$
K^a_{\v,q}(x/q)=K^a_{\v,q}(qx)=0,
$$
and then $K^a_{\v,q}(x)=0,\quad\forall x\in\R$, but this is impossible. This proves that $K^a_{\v,q}(x)>0,\quad\forall x\in\R.$
\end{proo}

\begin{coro}\label{c2}
Let $h$ a given function defined on $\R$ then we have
$$
\left[1-\frac{\Delta_{q,\v}}{a^2}\right]h(x) =-\frac{q^{2\v-1}(1-q)^2}{a^2x^{2\v+1}K_{\v,q}^a(x)%
}\Lambda_q^{-1}D_q\left[x^{2\v+1}K_{\v,q}^a(x)K_{\v,q}^a(qx)D_q\left[\frac{h(x)}{K_{\v,q}^a(x)}
\right]\right],
$$
\end{coro}

\begin{proo}
As we have
$$
\left[1-\frac{\Delta_{q,\v}}{a^2}\right]K_{\v,q}^a=0.
$$
Then the proof is identical to that of Proposition \ref{p1}.
\end{proo}

\section{Composite kernels}

In this section we study the composite variation diminishing $*_q$-kernel
using the results proved in preceding section.

\begin{lemm}\label{lem}If $g_1$ and $g_2$ are variation diminishing $*_q$-kernels
then $g=g_1*_qg_2$ is also a variation diminishing $*_q$-kernel.
\end{lemm}

\begin{proo}Let $f\in\cb$, then by (\ref{e12}), the function $g_2*_qf\in\cb$ and therfore
$$
V[g*_qf]=V[g_1*_q(g_2*_qf)]\leq V[g_2*_qf]\leq V[f].
$$
which achieves the proof.
\end{proo}

\begin{theore}\label{t3}
Let  $0<a_1\leq a_2\leq\ldots$ where
$$
\sum_{k=1}^\I a_k^{-2}<\I.
$$
If $\ds E(t)=\prod_{k=1}^\I\left[1+\frac{t^2}{a_k^2}\right]$ then $1/E$ is the $q$-Bessel Fourier transform of a positive variation
diminishing $*_q$-kernel $G$, and we have
$$
\F(G)(x)=1/E(x),\quad\forall x\in\R.
$$
In particular the $q$-Gauss kernel $e^c_{\v,q}$ is a variation diminishing $*_q$-kernel.
\end{theore}

\begin{proo}
Let  $ G_n(t)=K_{\nu}^{a_1}*_q\ldots*_qK_{\nu}^{a_n}(t)$.  By Theorem \ref{t1} and (\ref{e5}) we have $G_n\in\ll$.
By (\ref{e6})
$$
\F(G_n)(t)=\prod_{i=1}^n\F(K_{\nu}^{a_i})(t)=\prod_{i=1}^n\left(1+\frac{x^2}{a_i^2}\right)^{-1}=\frac{1}{E_n(t)},
$$
where
$$
E_n(t)=\prod_{i=1}^n\left(1+\frac{x^2}{a_i^2}\right).
$$
By Theorem \ref{th} and Lemma \ref{lem} we see that $G_n$ is a variation diminishing $*_q$-kernel. On the other hand $\ds\frac{1}{E_n}\in\ll$ and therefore by the inversion formula  $\ds\F\left(\frac{1}{E_n}\right)=G_n$. By Lemma \ref{lem3} and the fact that $\F(G_n)(0)=1$ we have
\begin{equation}\label{e19}
G_n(x)\geq 0,\quad\forall x\in\R,\forall n\in\N^*.
\end{equation}
Now consider
$$
E(t)=\prod_{k=1}^\I\left[1+\frac{t^2}{a_k^2}\right].
$$
Since $1/E\in\ll$ we can introduce the following function
$$
G(x)=c_{q,\v}\int_0^\I\frac{1}{E(t)}j_\v(xt,q^2)t^{2\v+1}d_qt,
$$
which satisfies
$$
\F(G)=1/E.
$$
The function $G$ belong to $\ll$. Indeed $\ds\lim_{x\to 0}G(x)$ exists and we use  the same method to prove (\ref{e20}) . Also same arguments led to
$$
G_n\in\ll,\quad\forall n\in\N^*.
$$
On the other hand by  (\ref{e11})
$$
\|G-G_n\|_{q,\I}=\|\F(1/E)-\F(1/E_n)\|_{q,\I}\leq B_{q,\v}\|1/E-1/E_n\|_{q,1,\v}.
$$
Since $\ds\lim_{n\to\I}\|1/E-1/E_n\|_{q,1,\v}=0$ and by (\ref{e19}) we obtain
$$\lim_{n\to\I}\|G-G_n\|_{q,\I}=0\Rightarrow G(x)\geq0,\quad\forall x\in\R.
$$
Let $f\in\A$. From (\ref{e12}) we have
$$
\|G*_qf-G_n*_qf\|_{q,\I}\leq \|G-G_n\|_{q,\I}\|f\|_{q,1,\v}.
$$
Which implies
$$\lim_{n\to\I}\|G*_qf-G_n*_qf\|_{q,\I}=0,
$$
and thus for $n$ big enough we get
$$
V[G*_qf]=V[G_n*_qf]\leq V[f].
$$
Same technics as in Theorem \ref{th} leads to
$$
V[G*_qf]\leq V[f],\quad\forall f\in\cb.
$$
These means that $G$ is a variation diminishing $*_q$-kernel.\\

For the particular cases $a_k=cq^{-k}$ we have
$$
E(t)=\dfrac{1}{e(-c^2t^2,q^2)},\quad c> 0.
$$
By (\ref{e17}) we prove that the $q$-Gauss kernel $e^c_{\v,q}$ is a variation diminishing $*_q$-kernel.
\end{proo}

\begin{remar}
For $n\geq 2$ we have  $G_n=K_{\v,q}^{a_n}*_qG_{n-1}$. Corollary \ref{c1} implies that $G_n$ is a solution of the following $q$-difference equation:
$$
\left[1-\frac{\Delta_{q,\nu}}{a_n^2}\right]G_n(x)=G_{n-1}(x);\quad \forall x\in\R.
$$
Thus
$$
G_n(x)-G_{n-1}(x)=\frac{\Big[G_n(q^{-1}x)+q^{2\nu}G_n(qx)\Big]-\Big[1+q^{2\nu}\Big]G_n(x)}{(a_nx)^2},
$$
We will prove that $G_n-G_{n-1}\geq0$, we can show it by recurrence. In fact, if $n=2$  we have
$$
G_2(x)=K_{\v,q}^{a_2}*_qK_{\v,q}^{a_1}(x)\geq K_{\v,q}^{a_1}(x).
$$
Now if we assume that it is true to order $n-1$, we have
\begin{eqnarray*}
G_n(q^{-1}x)+q^{2\nu}G_n(qx)&=&K_{\v,q}^{a_n}*_qG_{n-1}(q^{-1}x)+q^{2\nu}K_{\v,q}^{a_n}*_qG_{n-1}(qx)\\
&=&K_{\v,q}^{a_n}*_q\big[G_{n-1}(q^{-1}x)+q^{2\nu}G_{n-1}(qx)\big]\\
&=&c_{q,\nu}\int_0^{\I}K_{\v,q}^{a_n}(y)\big[T_{q,q^{-1}x}^{\nu}G_{n-1}(y)+q^{2\nu}T_{q,qx}^{\nu}G_{n-1}(y)\big]y^{2\nu+1}d_qy\\
&=&c_{q,\nu}\int_0^{\I}K_{\v,q}^{a_n}(y)\big[T_{q,y}^{\nu}G_{n-1}(q^{-1}x)+q^{2\nu}T_{q,y}^{\nu}G_{n-1}(qx)\big]y^{2\nu+1}d_qy\\
&=&c_{q,\nu}\int_0^{\I}K_{\v,q}^{a_n}(y)T_{q,y}^{\nu}\big[G_{n-1}(q^{-1}x)+q^{2\nu}G_{n-1}(qx)\big]y^{2\nu+1}d_qy\\
&\geq& c_{q,\nu}\int_0^{\I}K_{\v,q}^{a_n}(y)T_{q,y}^{\nu}\big[(1+q^{2\nu})G_{n-1}(x)\big]y^{2\nu+1}d_qy\\
&\geq& (1+q^{2\nu})c_{q,\nu}\int_0^{\I}K_{\v,q}^{a_n}(y)T_{q,x}^{\nu}G_{n-1}(y)y^{2\nu+1}d_qy\\
&\geq& (1+q^{2\nu})\big[K_{\v,q}^{a_n}*_qG_{n-1}\big](x)=(1+q^{2\nu})G_{n}(x).
\end{eqnarray*}
Hence, $G_n(x)\geq G_{n-1}(x)$. On the other hand
$$
G(x)=\lim_{n\rightarrow\I}G_n(x).
$$
As $(G_n)_{n\geq0}$ is a positive and an increasing sequences then $G(x)\geq G_n(x); \quad\forall x\in\R$. Beppo-Levi's Theorem leads to
\begin{equation}\label{e13}
\lim_{n\to\I}\|G-G_n\|_{q,1,\v}=\int_0^{\I}\lim_{n\to\I}|G(t)-G_n(t)|t^{2\nu+1}d_qt=0.
\end{equation}
Given $f\in\cb$. From (\ref{e12}) we get
$$
\|G*_qf-G_n*_qf\|_{q,\I}\leq \|G-G_n\|_{q,1,\v}\|f\|_{q,\I}.
$$
By  (\ref{e13}) we see that
$$\lim_{n\to\I}\|G*_qf-G_n*_qf\|_{q,\I}=0,
$$
and thus for $n$ big enough we get
$$
V[G*_qf]=V[G_n*_qf]\leq V[f],
$$
which means that $G$ is a variation diminishing $*_q$-kernel.
\end{remar}

\section{Asymptotic expansion}

In this section we discuss the asymptotic expansion at infinity of the variation diminishing $*_q$-kernel .
\begin{propositio}\label{p2}
The $q$-Macdonald function $K^a_\v$ satisfies the following properties

a. for all $x\in\R$ we have
$$
\Lambda_q^{-1}D_qK^a_{\v,q}(x)=-\frac{q}{1-q}xK^a_{\v+1,q}(x).
$$

b. there exists a constante $d^a_\v>0$ such thta for all $x\in\R$  we have $$
x^{2(\v+1)}\Big[\frac{1}{1-q^{2\v+2}}K_{\v,q}(x)I^a_{\v+1,q}(x)+K^a_{\v+1,q}(x)I^a_{\v,q}(x)\Big]=d^a_\v.
$$

c. if $\v>0$ then we have
\begin{equation}\label{e14}
\lim_{x\to 0^+}x^{2\v}K^a_{\v,q}(x)=d^a_{\v-1}.
\end{equation}
\end{propositio}

\begin{proo}

a. using the fact that \cite[p.60]{D2}
$$
D_q\Big[j_\v(.,q^2)\Big](t)=-\frac{q^2}{(1-q)(1-q^{2\v+2})}tj_{\v+1}(qt,q^2),
$$
with $c_{q,\nu}=(1-q^{2\nu+2})c_{q,\nu+1}$ we deduce that
\begin{eqnarray*}
\Lambda_q^{-1}D_qK^a_{\v,q}(x)&=&c_{q,\nu}\int_0^{+\infty}\left[1+\frac{t^2}{a^2}\right]^{-1}\Lambda_q^{-1}D_qj_{\nu}(tx,q^2)t^{2\nu+1}d_qt
\\&=&-\frac{q}{(1-q)(1-q^{2\nu+2})}xc_{q,\nu}\int_0^{+\infty}\left[1+\frac{t^2}{a^2}\right]^{-1}j_{\nu+1}(tx,q^2)t^{2\nu+3}d_qt
\\&=&-\frac{q}{1-q}xK^a_{\nu+1,q}(x).
\end{eqnarray*}
b. By  (\ref{e4}) and (\ref{e18}) and Remark \ref{re} we have:
$$
D_q\Big[y\mapsto y^{2\nu+1}w_y(I^a_{\nu,q},K^a_{\v,q})\Big]=\Big[\Delta_{q,\nu}I^a_{\nu}(x)K^a_{\v,q}(x)-I^a_{\nu,q}(x)\Delta_{q,\nu}K^a_{\v,q}(x)\Big]x^{2\nu+1}=0.
$$
c. By (\ref{e21})we have $\ds\lim_{x\to 0}x^{2\nu+2}K^a_{\v,q}(x)=0$. From (b.) and the fact that $\ds\lim_{x\to 0}I^a_{\nu,q}(x)=1$ we get
$$
\lim_{x\to 0}x^{2(\nu+1)}K^a_{\nu+1,q}(x)=d^a_\v,\quad \v>0,
$$
which gives the result.
\end{proo}

\begin{theore}
Let $\v>0$ and  $G$ be a  variation diminishing $*_q$-kernel such that $\ds\lim_{x\to 0}G(x)$ exists. There exists $a\in\R$ such that at infinity
$$
G(x)={\mathcal O}\Big[K^a_{\v,q}(x)\Big].
$$
As a consequence
\begin{equation}\label{e25}
G(q^n)={\mathcal O}\left(q^{n^2}\right),\quad\forall n\in\Z.
\end{equation}
\end{theore}

\begin{proo}By Lemma \ref{lem3} we can suppose without loss of generality that $G$ is positive on $\R$. We will prove that $G/K^a_\v$ does not have a local minimum in $\R$.\\

First note that the $q$-Gauss kernel  $e^c_{\v,q}>0$ on $\R$ and
$$
e^c_{\v,q}(qx)=\left(1+\frac{x^2}{q^{2\v}c}\right)e^c_{\v,q}(x),\quad e^c_{\v,q}(x/q)=\left(1+\frac{x^2}{q^{2\v+2}c}\right)^{-1}e^c_{\v,q}(x),
$$
which implies
$$
\left[1-\frac{\Delta_{q,\v}}{a^2}\right]e^c_{\v,q}(x)=\frac{1}{(ax)^2}\left[(ax)^2+\left(1+q^{2\v+2}\right)
-q^{2\v}\left(1+\frac{x^2}{q^{2\v}c}\right)-\left(1+\frac{x^2}{q^{2\v+2}c}\right)^{-1}\right]e^c_{\v,q}(x).
$$
Then $\ds\left[1-\frac{\Delta_{q,\v}}{a^2}\right]e^c_{\v,q}$ has at most one change of sign in $\R$. It follows since $G$ is a
variation diminishing $*_q$-kernel that
$$
G*_q\left[1-\frac{\Delta_{q,\v}}{a^2}\right]e^c_{\v,q}
$$
has at most one change of sign. On the other hand by (\ref{e23}) and (\ref{e6})
$$
\F\left[G*_q\left[1-\frac{\Delta_{q,\v}}{a^2}\right]e^c_{\v,q}\right]
=\F\left[\left[1-\frac{\Delta_{q,\v}}{a^2}\right]G*_qe^c_{\v,q}\right].
$$
By the inversion formula (\ref{e1}) we deduce that
$$
\left[1-\frac{\Delta_{q,\v}}{a^2}\right]G*_qe^c_{\v,q}=G*_q\left[1-\frac{\Delta_{q,\v}}{a^2}\right]e^c_{\v,q}.
$$
Let $c\to 0$, then $\ds\left[1-\frac{\Delta_{q,\v}}{a^2}\right]G$ has at most one change of sign. Furthermore by Corollary \ref{c2}
$$
\left[1-\frac{\Delta_{q,\v}}{a^2}\right]G(x)=-\frac{q^{2\v}}{a^2x^{2\v+1}K^a_{\v,q}(x)
}\Lambda_q^{-1}D_qH(x),
$$
where
$$
H(x)=x^{2\v+1}K^a_{\v,q}(x)K^a_{\v,q}(qx)D_q\left[\frac{G(x)}{K^a_{\v,q}(x)}\right]=x^{2\v+1}\Big[K^a_{\v,q}(x)D_qG(x)-G(x)D_qK^a_{\v,q}(x)\Big].
$$
We have
\begin{equation}\label{e15}
V\left[\left[1-\frac{\Delta_{q,\v}}{a^2}\right]G\right]=V\Big[D_qH\Big].
\end{equation}
Now if $G/K^a_{\v,q}$  hase a local minimum at $\delta\in\R$ then $H$
would have some negative values to the left of $\delta$ and some positive values to the right of $\delta$ and then $D_qH$ would have some
positive values near $\delta$.\\

The fact that $K^a_{\v,q}$ and $G$ belong to $\ll$ leads to
$$
\lim_{x\to\I}x^{2\v+2}K^a_{\v,q}(x)=\lim_{x\to\I}x^{2\v+2}G(x)=0\Rightarrow \lim_{x\to\I}H(x)=0.
$$
Then $D_qH$ would have some negative values near infinity. On the other hand by (\ref{e14})
$$
\lim_{x\to 0}x^{2\v+1}K^a_{\v,q}(x)D_qG(x)=\frac{d^a_{\v-1}}{1-q}\lim_{x\to 0}\Big[G(x)-G(qx)\Big]=0,
$$
and
$$
\lim_{x\to 0}x^{2\v+1}G(x)D_qK^a_{\v,q}(x)=\frac{G(0)}{1-q}\lim_{x\to 0}x^{2\v}\Big[K^a_{\v,q}(x)-K^a_{\v,q}(qx)\Big]=\frac{G(0)}{1-q}(1-q^{-2\v})d^a_{\v-1}
$$
which gives
$$
\lim_{x\to 0}H(x)\geq0.
$$
Since  $H$ hase some  negative values to the left of $\delta$ then $D_qH$ would have some
negative values near $0$. So,  $D_qH$ hase at least twos change of sign. But by (\ref{e15}) this is impossible.\\

Now $G/K^a_{\v,q}$ does not have a local minimum in $\R$ then  it's either non-decreasing on $\R$ or non-increasing. In the second case the theorem holds true. Suppose $G/K^a_{\v,q}$ is non-decreasing on $\R$ for all $a\in\R$. As we have
$$
K^a_{\v,q}(x)=a^{2\v+2}K_{\v,q}(ax),
$$
by (\ref{e14}) we obtain for a given $x\in\R$
$$
\lim_{a\to 0}\frac{x^{2\v}}{a^2}K^a_{\v,q}(x)=\lim_{a\to 0}(ax)^{2\v}K_{\v,q}(ax)=d^1_{\v-1}.
$$
then
$$
\lim_{a\to 0}a^2d^1_{\v-1}\frac{G(x)}{K^a_{\v,q}(x)}=x^{2\v}G(x).
$$
So $x^{2\v}G(x)$ is non-decreasing on $\R$ but this is impossible because $G\in\ll$.\\

To prove (\ref{e25}) we use Theorem 3 in \cite{D3} and (\ref{e4}) we obtain
$$
K^a_{\v,q}(q^n)={\mathcal O}\left(q^{n^2}\right),\quad\forall n\in\Z,
$$
which gives the result.
\end{proo}


\begin{thebibliography}{2}

\bibitem{D} L. Dhaouadi, On the $q$-Bessel Fourier transform,
Bulletin of Mathematical Analysis and Applications,
Volume 5, Issue 2, Article 42-60 (2013).

\bibitem{D1} L. Dhaouadi, A. Fitouhi and J. El Kamel,
Inequalities in $q$-Fourier Analysis,
Journal of Inequalities in Pure and Applied Mathematics,
Volume 7, Issue 5, Article 171 (2006).

\bibitem{D2} L. Dhaouadi, W. Binous and A. Fitouhi,
Paley-Wiener theorem for the $q$-Bessel transform and associated $q$-sampling formula,
Expo. Math,
Volume 27, Number 1, Article 55-72 (2009).


\bibitem{D3}L. Dhaouadi,
On q-Bessel Fourier analysis method for classical moment problem,
Boll. Unione Mat. Ital,
DOI: 10.1007/s40574-016-0115-8.


\bibitem{G} G. Gasper and M. Rahman,
Basic hypergeometric series,
Encycopedia of mathematics and its applications. Cambridge university press,
Volume  35 (1990).

\bibitem{H} I. I. Hirschman,  Jr,
Variation diminishing Hankel transforms,
J. Analyse. Math,
Volume 8, Article 307-336 (1961).

\bibitem{J} F. H. Jackson,
On a $q$-Definite Integrals,
Quarterly Journal of Pure and Application Mathematics,
Volume 41, Article 193-203 (1910).

\bibitem{K} T. H. Koornwinder and R. F. Swarttouw,
On $q$-Analogues of the Hankel and Fourier Transform,
Trans. A. M. S,
 Volume 333, Article 445-461 (1992).

\bibitem{K2} S. Karlin,
Total Positivity, Stanford University Press,
Stanford, Calif., 1968.

\bibitem{S} R. F. Swarttouw,
The Hahn-Exton $q$-Bessel functions,
PhD Thesis,
The Technical University of Delft (1992).

\bibitem{S2} I. J. Schoenberg, On P\`{o}lya frenquency function II,
Acta. Sci. Math. Szegzd,
Volume 12,  Article 97-106 (1950).


\bibitem{S1} I. J. Schoenberg,
On P\`{o}lya frenquency function I,
J. Analyse. Math,
Volume 1,  Article 331-374 (1951).



\bibitem{W}G. N. Watson,
A Treatise on the Theory of Bessel Functions,
Cambridge at the University Press (1922).

\end{thebibliography}
\end{document}